\documentclass[11pt]{amsart}
\usepackage{amsmath, amsfonts, eucal, amssymb, amsbsy}

\begin{document}
\newtheorem{thm}{Theorem}[section]
\newtheorem{lmm}[thm]{Lemma}
\newtheorem{cor}[thm]{Corollary}
\newtheorem{prop}[thm]{Proposition}
\newtheorem{defn}[thm]{Definition}
\newcommand{\argmax}{\operatorname{argmax}}
\newcommand{\argmin}{\operatorname{argmin}}
\newcommand{\bbb}{\mathbf{B}}
\newcommand{\bbr}{\mathbf{R}}
\newcommand{\bbw}{\mathbf{W}}
\newcommand{\bbx}{\mathbf{X}}
\newcommand{\bbxa}{\mathbf{X}^*}
\newcommand{\bbxb}{\mathbf{X}^{**}}
\newcommand{\bbxp}{\mathbf{X}^\prime}
\newcommand{\bbxpp}{\mathbf{X}^{\prime\prime}}
\newcommand{\bbxt}{\tilde{\mathbf{X}}}
\newcommand{\bby}{\mathbf{Y}}
\newcommand{\bbz}{\mathbf{Z}}
\newcommand{\bbzt}{\tilde{\mathbf{Z}}}
\newcommand{\bp}{b^\prime}
\newcommand{\bx}{\mathbf{x}}
\newcommand{\by}{\mathbf{y}}
\newcommand{\cc}{\mathbb{C}}
\newcommand{\cov}{\mathrm{Cov}}
\newcommand{\dd}{\mathcal{D}}
\newcommand{\ee}{\mathbb{E}}
\newcommand{\fp}{f^\prime}
\newcommand{\fpp}{f^{\prime\prime}}
\newcommand{\fppp}{f^{\prime\prime\prime}}
\newcommand{\ii}{\mathbb{I}}
\newcommand{\gp}{g^\prime}
\newcommand{\gpp}{g^{\prime\prime}}
\newcommand{\gppp}{g^{\prime\prime\prime}}
\newcommand{\ma}{\mathcal{A}}
\newcommand{\mf}{\mathcal{F}}
\newcommand{\mi}{\mathcal{I}}
\newcommand{\ml}{\mathcal{L}}
\newcommand{\mx}{\mathcal{X}}
\newcommand{\mxp}{\mathcal{X}^\prime}
\newcommand{\mxpp}{\mathcal{X}^{\prime\prime}}
\newcommand{\my}{\mathcal{Y}}
\newcommand{\myp}{\mathcal{Y}^\prime}
\newcommand{\mypp}{\mathcal{Y}^{\prime\prime}}
\newcommand{\pp}{\mathbb{P}}
\newcommand{\ppr}{p^\prime}
\newcommand{\pppr}{p^{\prime\prime}}
\newcommand{\ra}{\rightarrow}
\newcommand{\rr}{\mathbb{R}}
\newcommand{\sss}{\sigma^\prime}
\newcommand{\st}{\sqrt{t}}
\newcommand{\sst}{\sqrt{1-t}}
\newcommand{\tr}{\operatorname{Tr}}
\newcommand{\uu}{\mathcal{U}}
\newcommand{\var}{\mathrm{Var}}
\newcommand{\ve}{\varepsilon}
\newcommand{\vp}{\varphi^\prime}
\newcommand{\vpp}{\varphi^{\prime\prime}}
\newcommand{\ww}{W^\prime}
\newcommand{\xp}{X^\prime}
\newcommand{\xpp}{X^{\prime\prime}}
\newcommand{\xt}{\tilde{X}}
\newcommand{\xx}{\mathcal{X}}
\newcommand{\yp}{Y^\prime}
\newcommand{\ypp}{Y^{\prime\prime}}
\newcommand{\zt}{\tilde{Z}}
 
\newcommand{\fpar}[2]{\frac{\partial #1}{\partial #2}}
\newcommand{\spar}[2]{\frac{\partial^2 #1}{\partial #2^2}}
\newcommand{\mpar}[3]{\frac{\partial^2 #1}{\partial #2 \partial #3}}
\newcommand{\tpar}[2]{\frac{\partial^3 #1}{\partial #2^3}}

\title{Stein's method for concentration inequalities}
\author{Sourav Chatterjee}
\address{\newline367 Evans Hall \#3860\newline
Department of Statistics\newline
University of California\newline
Berkeley, CA 94720-3860\newline
{\it E-mail: \tt sourav@stat.berkeley.edu}\newline 
{\it Homepage: \tt http://www.stat.berkeley.edu/$\sim$sourav}
}
\maketitle

\begin{abstract}
We introduce a version of Stein's method for proving  concentration and moment inequalities in problems with dependence. Simple illustrative examples from combinatorics, physics, and mathematical statistics are provided.
\end{abstract}

\section{Introduction and results}\label{intro}
Stein's method was introduced by Charles Stein \cite{stein72} in the context of normal approximation for sums of dependent random variables. Stein's version of his method, best known as the ``method of exchangeable pairs'', attained maturity in his later work \cite{stein86}. A reasonably large literature has developed around the subject, but it has almost exclusively developed as a method of proving distributional convergence with error bounds. Stein's attempts at getting large deviations in \cite{stein86} did not, unfortunately, prove fruitful. Some progress for sums of dependent random variables was made by Rai\v{c} \cite{raic04}. A~general version of Stein's method for concentration inequalities was introduced for the first time in the Ph.D.\ thesis \cite{chatterjee05} of the present author. The purpose of this paper is to explain the theory developed in \cite{chatterjee05} via examples. Another application is in \cite{chatterjee05a}.

This section is organized as follows: First, we give three examples, followed by the main abstract theorem; finally, towards the end of the section, we present very condensed overviews of Stein's method, concentration of measure, and the related literature. Proofs are in section \ref{proofs}.

\subsection{A generalized matching problem}\label{perm}
Let $\{a_{ij}\}$ be an $n\times n$ array of real numbers. Let
$\pi$ be chosen uniformly at random from the set of all permutations of $\{1,\ldots,n\}$, and let $X=\sum_{i=1}^n a_{i\pi(i)}$. This class of random variables was first studied by Hoeffding \cite{hoeffding51}, who proved that they are approximately normally distributed under certain conditions. It is easy to see that various well-studied functions of random permutations, like the number of fixed points, the sum of a random sample picked without replacement from a finite population, and the function $\sum_i|i-\pi(i)|$ (known as Spearman's footrule~\cite{diaconisgraham77}), are all instances of Hoeffding's statistic. 

Hoeffding's statistic has a long history of association with Stein's method. In fact, in an unpublished work Stein introduced his method to treat the normal approximation problem for this object. Bolthausen \cite{bolthausen84} used Stein's method to give a Berry-Esseen bound. Bolthausen and G\"otze \cite{bolthausengotze93} gave  multivariate central limit theorems under a further generalized setup. However, we have not seen large deviations or concentration bounds using any method.

Our version of Stein's method enables us to easily derive the following nice tail bound.
\begin{prop}\label{hoeff3}
Let $\{a_{ij}\}_{1\le i,j\le n}$ be a collection of numbers from $[0,1]$. Let $X = \sum_{i=1}^n a_{i\pi(i)}$, where $\pi$ is drawn from the uniform distribution over the set of all permutations of $\{1,\ldots,n\}$. Then
\[
\pp\{|X - \ee(X)| \ge t\} \le
2\exp\biggl(-\frac{t^2}{4\ee(X) + 2t}\biggr)
\]
for any $t\ge 0$.
\end{prop}
\noindent Note that the bound does not have an explicit dependence on $n$. Note also the automatic transition from Poissonian to gaussian tails as $\ee(X)$ becomes large (when $\ee(X)$ is small the bound is like $\exp(-Ct)$, whereas when $\ee(X)$ is large, it is essentially a gaussian tail with standard deviation $\sqrt{\ee(X)}$.). These two properties characterize it as a so-called ``Bernstein type inequality'', named after the classical Bernstein inequality (see \cite{shorackwellner86}, page 855) for sums of bounded independent random variables.

The classical result of Maurey \cite{maurey79} can only imply the weaker inequality $P(X> \ee(X) +t) \le e^{-t^2/4n}$. However, it is possible to derive a Bernstein bound similar to Proposition \ref{hoeff3} (albeit with a significantly worse constant in the exponent) using Michel Talagrand's deep theorem about concentration of random permutations (Theorem 5.1 in Section~5 of~\cite{talagrand95}; see also McDiarmid~\cite{mcdiarmid02} and Luczak \& McDiarmid~\cite{luczakmcdiarmid03}). 

For a concrete application, let $X$ be the number of fixed points of a random permutation $\pi$. Then $X =\sum_{i=1}^n a_{i\pi(i)}$, where $a_{ij} = \ii_{\{i=j\}}$. Since $\ee(X)=1$, Proposition \ref{hoeff3} gives $\pp\{|X-1|\ge t\} \le 2\exp(-t^2/(4+2t))$. Of course, we do not expect this to be the best possible bound in this very well-understood  problem; this is just meant to be an illustration. In fact, the exact distribution of the the number of fixed points is known (see Feller~\cite{feller1}, section IV.4), which gives a tail bound like $\exp(-Ct\log t)$.

Finally, we also have a ``Burkholder-Davis-Gundy'' type inequality for Hoeffding's statistic which does not require a bound on the $a_{ij}$'s.
\begin{prop}\label{hoeff2}
Let $\{a_{ij}\}_{1\le i,j\le n}$ be an arbitrary collection of real numbers. Let $\pi$ be a uniform random permutation, and let $X = \sum_{i=1}^n a_{i\pi(i)}$. Define
\[
\Delta = \frac{1}{4n}\sum_{i,j} (a_{i\pi(i)} +
a_{j\pi(j)}-a_{i\pi(j)}-a_{j\pi(i)})^2.
\]
Then for every positive integer $k$, we have $\ee(X-\ee(X))^{2k} \le (2k-1)^k \ee \Delta^k$.
\end{prop}
\noindent For a general exposition about the famous  Burkholder-Davis-Gundy martingale inequalities we refer to the article by Burkholder \cite{burkholder73}.

\subsection{Magnetization in the Curie-Weiss model}\label{curieweiss}
Fix any $\beta \ge 0$, $h\in \rr$, and consider the probability mass function (the Gibbs measure) on $\{-1,1\}^n$ given by
\begin{equation}\label{gibbs}
\pp(\{\sigma\}) := Z^{-1}\exp\biggl(\frac{\beta}{n}\sum_{i<j} \sigma_i\sigma_j +\beta h \sum_i \sigma_i\biggr),
\end{equation}
where $\sigma = (\sigma_1,\ldots,\sigma_n)$ is a typical element of $\{-1,1\}^n$ and $Z$ is the normalizing constant (depends on $\beta$ and $h$). This is known as the `Curie-Weiss model of ferromagnetic interaction' at inverse temperature $\beta$ and external field $h$. The $\sigma_i$'s stand for the spins of $n$ particles, each having a spin of $+1$ or $-1$. The ferromagnetic interaction between the particles is captured in a very simplistic manner by the first term in the hamiltonian. 

The {\it magnetization} of the system, as a function of the configuration $\sigma$, is defined as
$m(\sigma) := \frac{1}{n}\sum_{i=1}^n \sigma_i$. 
If $n$ is large and $\sigma$ is drawn from the Gibbs measure, then the magnetization satisfies 
\begin{equation}\label{meq}
m(\sigma) \approx \tanh(\beta m(\sigma) + \beta h).
\end{equation}
with high probability. 
The equation has a unique root for small values of $\beta$ and multiple solutions for $\beta$ above a critical value. In the physics parlance, this is described by saying that the Curie-Weiss model exhibits ``spontaneous magnetization'' at low temperatures. For a formal discussion with rigorous proofs, we refer to Ellis \cite{ellis85}, section IV.4.

The following proposition formalizes \eqref{meq} with finite sample tail bounds. 
\begin{prop}\label{curie}
Suppose $\sigma$ is drawn from the Gibbs measure \eqref{gibbs}. Then, for any $\beta \ge 0$, $h\in \rr$, $n\ge 1$, and $t\ge 0$, the magnetization $m := \frac{1}{n}\sum_i \sigma_i$ satisfies 
\[
\pp\biggl\{\bigl|m -  \tanh(\beta m + \beta h)\bigr| \ge \frac{\beta}{n} + \frac{t}{\sqrt{n}}\biggr\} \le 2\exp\biggl(-\frac{t^2}{4(1+\beta)}\biggr).
\]
\end{prop}
\noindent Although the Curie-Weiss model is a simple model of ferromagnetic interaction, we haven't encountered any result in the literature which gives an explicit bound like the above. In particular, the result shows concentration of $m(\sigma)$ around the {\it set} of roots of $x=\tanh(\beta x +\beta h)$, and not just its  mean.

However, concentration inequalities for Gibbs measures without explicit constants under various mixing conditions have been obtained before. For a history of the literature and some significant recent progress, we refer to Chazottes et.~al.~\cite{chazottes06}.

\subsection{Least squares estimation in the Ising model}\label{least}
The Ising model is another model of ferromagnetic interaction. Given an undirected graph $G = (V,E)$ on the vertex set $V = \{1,\ldots,n\}$, the Ising model without external field assigns the following probability density on $\{-1,1\}^n$:
\begin{equation}\label{ising}
\pp(\{\sigma\}) = Z(\beta)^{-1} \exp\biggl(\beta \sum_{\{i,j\}\in E } \sigma_i\sigma_j\biggr).
\end{equation}
Here, as before, $\beta$ is the inverse temperature and $Z(\beta)$ is the normalizing constant. 
A natural statistical problem in this model is the following: How to make inference about $\beta$ when your data is a single configuration generated from the Gibbs measure? 

The classical maximum likelihood approach for this problem was first considered by Pickard \cite{pickard87}. Iterative methods for computing the maximum likelihood estimator (e.g.\ Geyer \& Thompson \cite{geyerthompson92}, Jerrum \& Sinclair \cite{jerrumsinclair93}) are widely used nowadays. The Jerrum-Sinclair algorithm for computing the normalizing constant in the Ising model provably converges in polynomial time. However, it is not so clear whether the MLE is a good estimator at all, particularly at critical temperatures.

Here we investigate a method of estimating $\beta$ by minimizing an explicit sum-of-squares. First, let $\sigma$ be drawn from the Gibbs measure \eqref{ising} on $\{-1,1\}^n$, and for each $i$, let
\[
m_i := \sum_{j:\{i,j\}\in E} \sigma_j.
\]
For each $u \ge 0$, let
\begin{equation}\label{stheta}
S(u) := \frac{1}{n}\sum_{i=1}^n \bigl(\sigma_i - \tanh(u m_i))^2.
\end{equation}
The `least-squares estimate' of $\beta$ is defined to be
\[
\hat{\beta}_{LS} := \argmin_{u\ge 0} S(u). 
\]
Note that it is practically very easy to compute $\hat{\beta}_{LS}$, because $S$ is a smooth function of a single variable.

The least-squares technique is well-known and commonly used in the analysis of gaussian Markov random field (GMRF) models (probably originating from Besag \cite{besag75}), but rigorous results are scarce. 

Proposition \ref{isingprop} (stated below) shows that the random function $S$ indeed attains an approximate global minimum near $\beta$. In fact, it gives
\[
\ee|S(\beta) - \min_{u\ge 0} S(u)| = O\biggl(\sqrt{\frac{r\log n}{n}}\biggr),
\]
where $r$ is the maximum degree of the dependency graph $G$ (recall that the degree of a vertex is the number of neighbors of that vertex, and the maximum degree of a graph is the maximum vertex degree).
\begin{prop}\label{isingprop}
Let $r$ be the maximum degree of the dependency graph $G$ in the Ising model \eqref{ising}, and let $S(u)$ be defined as in \eqref{stheta}. Take any $t\ge 0$ and let
\[
\varepsilon = \sqrt{\frac{r(\log n + t)}{n}}.
\]
Then we have
\[
\pp\{S(\beta) \ge \min_{u\ge 0} S(u) + C\varepsilon\} \le \exp(-Kt^2),
\]
where $C$ and $K$ are numerical constants.
\end{prop}
\noindent Although it is unclear whether Proposition \ref{isingprop} is useful from a statistical point of view, it seems to be interesting as a mathematical result. For instance, observe that the conclusion is valid at any temperature. This is quite remarkable, since the low temperature phase in the Ising model is notoriously intractable for most graphs.

Here we should also mention that the technique can be easily applied to the Ising model with an external field, but we prefer to restrict ourselves to the problem of estimating a single parameter (the temperature) for the sake of clarity.

\subsection{The abstract result}\label{theory}
The following theorem encapsulates the concentration and moment inequalities used to work out all the examples in this paper.
\begin{thm}\label{conc}
Let $\xx$ be a separable metric space and suppose $(X,\xp)$ is an exchangeable pair of $\xx$-valued random variables. Suppose $f:\xx\ra \rr$ and $F:\xx\times \xx\ra \rr$ are square-integrable functions such that $F$ is antisymmetric (i.e.\ $F(X,\xp)=-F(\xp,X)$ a.s.), and $\ee(F(X,\xp)\mid X) = f(X)$ a.s. Let 
\[
\Delta(X) := \frac{1}{2}\ee\bigl(|(f(X)-f(\xp))F(X,\xp)|\,\bigl|\, X\bigr).
\]
Then $\ee(f(X)) = 0$, and the following concentration results hold for $f(X)$:
\begin{enumerate}
\item[$(i)$] If $\ee(\Delta(X)) <\infty$, then $\var(f(X)) = \frac{1}{2}\ee((f(X)-f(\xp))F(X,\xp))$.
\item[$(ii)$] Assume that $\ee(e^{\theta f(X)} |F(X,\xp)|) < \infty$ for all $\theta$. If there exists nonnegative constants $B$ and $C$ such that $\Delta(X) \le B f(X) + C$ almost surely, then for any $t\ge 0$, 
\[
\pp\{f(X) \ge t\}\le \exp\biggl(-\frac{t^2}{2C + 2Bt}\biggr) \ \ \text{and} \ \ \pp\{f(X) \le -t\}\le \exp\biggl(-\frac{t^2}{2C}\biggr).
\]
\item[$(iii)$] For any positive integer $k$, we have the following exchangeable pairs version of the Burkholder-Davis-Gundy inequality:
\[
\ee(f(X)^{2k}) \le (2k-1)^k \ee(\Delta(X)^k).
\] 
\end{enumerate}
\end{thm}
\noindent To see how the exchangeable pairs are constructed and the theorem is applied in our examples, one has to look at the proofs in section \ref{proofs}. However, for a quick illustration, we will now work out the inequalities for sums of independent random variables, taking care to spell out details.

\subsection{Simplest example} Let $X = \sum_{i=1}^n Y_i$, where $Y_i$'s are independent square integrable random variables. Let $\mu_i =\ee(Y_i)$ and $\sigma_i^2 = \var(Y_i)$. An exchangeable pair is created by choosing a coordinate $I$ uniformly at random from $\{1,\ldots,n\}$, and defining
\[
\xp = \sum_{j\ne I} Y_j + Y^\prime_I,
\]
where $Y^\prime_1,\ldots,Y^\prime_n$ are independent copies of $Y_1,\ldots,Y_n$. Let 
\[
F(x,y) = n(x-y).
\]
Then
\begin{align*}
\ee(F(X,\xp)\mid Y_1,\ldots,Y_n) &= \frac{1}{n}\sum_{i=1}^n \ee(n(Y_i - Y^\prime_i)\mid Y_1,\ldots,Y_n) \\
&= \sum_{i=1}^n (Y_i - \mu_i) = X - \ee(X).
\end{align*}
Since the right hand side depends only on $X$, we have 
\[
f(X) = \ee(F(X,\xp)\mid X) = X - \ee(X).
\]
Thus, from part ($i$) of Theorem \ref{conc} we get the elementary identity 
\[
\var(X) = \frac{1}{2}\sum_{i=1}^n \ee(Y_i - Y^\prime_i)^2 = \sum_{i=1}^n \sigma_i^2.
\]
Now note that
\begin{align*}
\Delta(X) &= \frac{n}{2}\ee((X-\xp)^2 \mid X) \\
&= \frac{1}{2}\sum_{i=1}^n\ee((Y_i-\yp_i)^2\mid X).
\end{align*}
If $c_1,\ldots,c_n$ are constants such that $|Y_i-\mu_i|\le c_i$ a.s. for each $i$, then 
\begin{align*}
\ee((Y_i-\yp_i)^2 \mid X) &= \ee((Y_i-\mu_i)^2\mid X) + \ee((\yp_i-\mu_i)^2)\\
&\le c_i^2 + \sigma_i^2.
\end{align*}
Part ($ii$) of Theorem \ref{conc} now implies that
\[
\pp\{|X- \ee(X)|\ge t\} \le 2 \exp\biggl(-\frac{t^2}{\sum_{i=1}^n (c_i^2 + \sigma_i^2)}\biggr).
\]
This is similar to (but not exactly the same as) the classical Hoeffding inequality \cite{hoeffding63} for sums of bounded random variables.

Now suppose that $0\le Y_i \le 1$ a.s.\ for each $i$. If the $\mu_i$'s are very small, then the Hoeffding bound is wasteful. A more careful analysis gives a better result, as follows. First, note that
\begin{align*}
\Delta(X) &= \frac{1}{2} \sum_{i=1}^n \ee((Y_i-\yp_i)^2\mid X)\\
&= \frac{1}{2}\sum_{i=1}^n (\ee Y_i^2 - 2\mu_i \ee(Y_i\mid X) + \ee(Y_i^2\mid X)).
\end{align*}
Using the assumption that $0\le Y_i\le 1$, we get
\[
\Delta(X) \le \frac{1}{2}\sum_{i=1}^n (\ee(Y_i) + \ee(Y_i\mid X)) = \frac{1}{2}(\ee(X) + X) = \frac{1}{2}f(X) + \ee(X).
\]
Thus, we can take $B = 1/2$ and $C=\ee(X)$ in part ($ii$) of Theorem \ref{conc}, which gives
\[
\pp\{|X-\ee(X)| \ge t\} \le 2\exp\biggl(-\frac{t^2}{2\ee(X)+t}\biggr).
\]
Again, this is a version of the classical Bernstein inequality (see \cite{shorackwellner86}, page 855) for sums of independent random variables.

Finally observe that by part ($iii$) of Theorem \ref{conc} and an application of Jensen's inequality, we have for each positive integer $k$,
\begin{align*}
\ee(X^{2k}) &\le (2k-1)^k \ee\biggl(\frac{1}{2}\sum_{i=1}^n \ee((Y_i-\mu_i)^2 + (\yp_i-\mu_i)^2\mid X)\biggr)^k\\
&\le (2k-1)^k \ee\biggl(\sum_{i=1}^n (Y_i-\mu_i)^2\biggr)^k.
\end{align*}
This is exactly what the Burkholder-Davis-Gundy inequality \cite{burkholder73} would give us for sums of independent random variables (although in this case, it can be derived by easier methods).

In the remainder of this section, we give very short overviews of Stein's method and concentration of measure. 

\subsection{Stein's method} 
Suppose we want to show that a random variable $X$ taking value in some space $\xx$ has approximately the same distribution as some other random variable $Z$. The classical version of Stein's method \cite{stein72, stein86} involves four steps:
\begin{enumerate}
\item Identify a ``characterizing operator'' $T$ for $Z$, which has the defining property that for any function $g$ belonging to a fixed large class of functions, $\ee Tg (Z) = 0$. For instance, if $\xx = \rr$ and $Z$ is a standard gaussian random variable, then $Tg(x) := g^\prime(x)-xg(x)$ is a characterizing operator, acting on all locally absolutely continuous $g$ with subexponential growth at infinity.  
\item Construct a random variable $\xp$ such that $(X,\xp)$ is an exchangeable~pair.
\item Find an operator $\alpha$ such that for any suitable $h:\xx \ra \rr$, $\alpha h$ is an antisymmetric function (i.e.\ $\alpha h(x,y)\equiv -\alpha h(y,x)$) and
\[
|\ee(\alpha h(X,\xp)|X = x) - Th(x)| \le \varepsilon_h,
\]
 where $\varepsilon_h$
is a small error depending only on $h$. 
\item Take a function $g$ and find $h$ such that $Th(x) = g(x) - \ee g(Z)$. By antisymmetry of $\alpha h$ and the exchangeability of $(X,\xp)$, it follows that $\ee(\alpha h(X,\xp))=0$. Combining with the previous step, we have the error bound $|\ee g(X) - \ee g(Z)|\le \varepsilon_h$. 
\end{enumerate}
There are other variants of Stein's method, most notably the generator method of Andrew Barbour \cite{barbour90}, the dependency graph approach introduced by Chen \cite{chen75} and Baldi and Rinott \cite{baldirinott89} and popularized by Arratia, Goldstein and Gordon \cite{agg90}, the size-biased coupling method of Barbour, Holst and Janson \cite{bhj92}, and the zero-biased coupling method due to Goldstein and Reinert \cite{goldsteinreinert97}. The recent applications to algebraic problems by Jason Fulman \cite{fulman04, fulman05}, and the quest for Berry-Esseen bounds by Rinott and Rotar \cite{rinottrotar97} and Shao and Su \cite{shaosu04} are also worthy of note. 

However, it is not our purpose here to go deeply into the regular versions of Stein's method. For further references and exposition, we refer to the recent monograph \cite{stein04}. For applications of the method of exchangeable pairs and other versions of Stein's method to Poisson approximation, one can look at the survey paper by Chatterjee, Diaconis \& Meckes \cite{cdm05}.

\subsection{Concentration inequalities}
The theory of concentration inequalities tries to answer the following question: Given a random variable $X$ taking value in some measure space $\xx$ (which is usually some high dimensional Euclidean space), and a measurable map $f:\xx \ra \rr$, what is a good explicit bound on $\pp\{|f(X)-\ee f(X)| \ge x\}$? Exact evaluation or accurate approximation is, of course, the central purpose of probability theory itself. In  situations where this is not possible, concentration inequalities aim to do the next best job by providing rapidly decaying tail bounds.

The literature on concentration inequalities is huge --- from the pioneering inequalities of Hoeffding \cite{hoeffding63} to the momentous work of Talagrand \cite{talagrand95} --- but most of it revolves around well-behaved functions of independent random variables. 
For a nearly complete account of the literature until the year 2001, we redirect the reader to the definitive resource in this subject --- the monograph~\cite{ledoux01} by Michel Ledoux. The methods of Kim and Vu \cite{kimvu04} and Boucheron, Lugosi, and Massart \cite{blm03} are significant recent developments.

The techniques developed in \cite{chatterjee05} (and partially presented here) have some basic similarities with the concentration results of Schmuckenschl\"ager \cite{schmuck98}, but go much beyond that in terms of applications. Other than that (and log-Sobolev inequalities, which are much harder to obtain anyway) there is very little --- even in the vast concentration literature  --- about the concentration of functions of dependent random variables, particularly in the discrete setting. We hope that our version of Stein's method will partially fill this~void.
\vskip.2in
\noindent {\bf Acknowledgments.} I am grateful to Persi Diaconis and Yuval Peres for many useful comments and suggestions. Thanks are also due to the two anonymous referees for pointing out several omissions and errors.

\section{Proofs}\label{proofs}
Before proving Theorem \ref{conc}, let us see how it is applied to work out the three examples described in section \ref{intro}.
\vskip.1in
\noindent {\bf Proof of Proposition \ref{hoeff3}.} Construct $\xp$ as follows: Choose $I,J$ uniformly and independently at random from $\{1,\ldots, n\}$. Let $\pi^\prime = \pi\circ (I,J)$, where $(I,J)$ denotes the transposition of $I$ and $J$. It can be easily verified that $(\pi,\pi^\prime)$ is an exchangeable pair. Hence if we let 
\[
\xp := \sum_{i=1}^n a_{i\pi^\prime(i)},
\]
then $(X,\xp)$ is also an exchangeable pair. Now note that
\begin{align*}
\frac{1}{2}\ee(n(X - \xp)|\pi) &= \frac{n}{2}\ee(a_{I\pi(I)} + a_{J\pi(J)} - a_{I\pi(J)} - a_{J\pi(I)}|\pi) \\
&= \frac{1}{n} \sum_{i,j} a_{i\pi(i)} - \frac{1}{n}\sum_{i,j} a_{i\pi(j)} \\
&= X-\ee(X).
\end{align*}
Thus, we can take $f(x)= x-\ee(X)$ and $F(x,y) = \frac{1}{2}n(x-y)$. Now note that since $0\le a_{ij}\le 1$ for
all $i$ and $j$, we have
\begin{align*}
&\frac{1}{2}\ee\bigl(|(f(X)-f(\xp))F(X,\xp)| \, \bigl| \,\pi\bigr) = \frac{n}{4}\ee((X-\xp)^2|\pi) \\
&=\frac{1}{4n}\sum_{i,j} (a_{i\pi(i)} +
a_{j\pi(j)}-a_{i\pi(j)}-a_{j\pi(i)})^2 \\
&\le \frac{1}{2n} \sum_{i,j}(a_{i\pi(i)} +
a_{j\pi(j)}+a_{i\pi(j)}+a_{j\pi(i)}) \\
&= X + \ee(X) = f(X) + 2\ee(X).
\end{align*}
Since the last quantity depends only on $X$ it follows that $\Delta(X)= f(X)+2\ee(X)$.
Applying part ($ii$) of Theorem \ref{conc} with $B = 1$ and $C= 2\ee(X)$ completes
the proof. \hfill $\Box$
\vskip.2in
\noindent {\bf Proof of Proposition \ref{hoeff2}.} Follows directly from part $(iii)$ of Theorem \ref{conc} and the computations done in the proof of Proposition \ref{hoeff3}.\hfill$\Box$
\vskip.2in
\noindent {\bf Proof of Proposition \ref{curie}.} Suppose $\sigma$ is drawn from the Gibbs distribution. We construct $\sss$ by taking a step in the Gibbs sampler as follows: Choose a coordinate $I$ uniformly at random, and replace the $I^{\mathrm{th}}$ coordinate of $\sigma$ by an element drawn from the conditional distribution of the $I^{\mathrm{th}}$ coordinate given the rest. 
It is well-known and easy to prove that $(\sigma,\sss)$ is an exchangeable pair. Let 
\[
F(\sigma, \sss) := \sum_{i=1}^n (\sigma_i - \sss_i).
\]
Now define
\[
m_i(\sigma) := \frac{1}{n}\sum_{j \le n, j\ne i} \sigma_j, \ \ i=1,\ldots,n.
\]
Since the Hamiltonian is a simple explicit function, the conditional distribution of the $i^{\mathrm{th}}$ coordinate given the rest is easy to obtain. An easy computation gives
$\ee(\sigma_i|\{\sigma_j, j\ne i\}) = \tanh(\beta m_i+\beta h)$. 
Thus, we have
\begin{align*}
f(\sigma) = \ee(F(\sigma,\sss)|\sigma) &= \frac{1}{n}\sum_{i=1}^n (\sigma_i -\ee(\sigma_i|\{\sigma_j, j\ne i\}))\\
&= m - \frac{1}{n}\sum_{i=1}^n \tanh(\beta m_i + \beta h).
\end{align*}
Now note that $|F(\sigma,\sss)| \le 2$, because $\sigma$ and $\sss$ differ at only one coordinate. Also, since the map $x\mapsto \tanh x$ is $1$-Lipschitz, we have 
\begin{align*}
|f(\sigma)-f(\sss)| &\le |m(\sigma)-m(\sss)| +\frac{\beta}{n}\sum_{i=1}^n |m_i(\sigma) - m_i(\sss)| \le \frac{2(1+\beta)}{n}.
\end{align*}
Thus, by part $(ii)$ of Theorem \ref{conc} we have
\[
\pp\biggl\{\biggl|m - \frac{1}{n}\sum_{i=1}^n \tanh(\beta m_i + \beta h)\biggr| \ge \frac{t}{\sqrt{n}}\biggr\} \le 2\exp\biggl(-\frac{t^2}{4(1+\beta)}\biggr).
\]
Finally note that for each $i$, by the Lipschitz nature of the $\tanh$ function, we get
\begin{align*}
&\biggl|\frac{1}{n}\sum_{i=1}^n \tanh(\beta m_i + \beta h) - \tanh(\beta m + \beta h) \biggr|\\
&\le \frac{1}{n}\sum_{i=1}^n|\tanh(\beta m_i+\beta h) - \tanh(\beta m+\beta h)| \\
&\le \frac{1}{n}\sum_{i=1}^n\beta|m_i-m| \le \frac{\beta}{n}.
\end{align*}
This completes the proof. \hfill $\Box$
\vskip.2in
\noindent{\bf Proof of Proposition \ref{isingprop}.}
As in the proof of Proposition \ref{curie}, we produce $\sss$ by taking a step in the Gibbs sampler: A coordinate $I$ is chosen uniformly at random, and $\sigma_I$ is replace by $\sss_I$ drawn from the conditional distribution of the $I^{\mathrm{th}}$ coordinate given $(\sigma_j)_{j\ne I}$. For each $i$, let 
\[
m_i = m_i(\sigma) := \sum_{j: \{i,j\}\in E} \sigma_j.
\]
Now fix $u\ge 0$ and define
\[
F(\sigma, \sss) := (\sigma_I- \sss_I) (\tanh(\beta m_I) - \tanh(u m_I)).
\]
Then $F(\sigma,\sss)= -F(\sss,\sigma)$ because $m_I(\sigma) = m_I(\sss)$ . Now let
\begin{align*}
&f(\sigma) := \ee(F(\sigma,\sss)\mid \sigma) \\
&= \frac{1}{n}\sum_{i=1}^n(\sigma_i - \tanh(\beta m_i)) (\tanh(\beta m_i) - \tanh(u m_i)).
\end{align*}
Now, if $r$ is the maximum degree of $G$, then at most $r+1$ terms in the sums defining $f(\sigma)$ and $f(\sss)$ are unequal, and they all lie in the interval $[-4,4]$. Thus, $|f(\sigma)-f(\sss)|\le 8(r+1)/n$. Also, evidently, $|F(\sigma, \sss)|\le 4$. Using all this information in part ($ii$) of Theorem \ref{conc}, we get
\[
\pp\{f(\sigma)\le -t\} \le \exp\biggl(-\frac{nt^2}{32(r+1)}\biggr).
\]
Now, a direct verification shows that
\begin{align*}
S(u) - S(\beta) &= \frac{1}{n}\sum_{i=1}^n (\tanh\beta m_i)-\tanh(u m_i))^2 + 2f(\sigma).
\end{align*}
Thus,
\begin{equation}\label{tail1}
\pp\{S(\beta) \ge S(u) + t\} \le \pp\{2f(\sigma) \le -t\} \le \exp\biggl(-\frac{nt^2}{128(r+1)}\biggr).
\end{equation}
Now note that for any $u,v\ge 0$, we have
\begin{align*}
&|S(u)-S(v)| \\
&\le \frac{1}{n}\sum_{i=1}^n |(2\sigma_i - \tanh(u m_i) - \tanh(vm_i))(\tanh(vm_i)-\tanh(um_i))|\\
&\le \frac{4}{n}\sum_{i=1}^n |\tanh(vm_i) - \tanh(um_i)| \le 4r|u-v|,
\end{align*}
since $|m_i(u-v)|\le r|u-v|$. 
Let $N = \lfloor\sqrt{nr\log n}\rfloor$, and let 
\[
u_k = k\sqrt{\frac{\log n}{nr}} \ \text{ for } \ k=1,2,\ldots,N.
\]
Then, if $u_{k-1} \le u\le u_k$, the above inequality gives
\[
|S(u)-S(u_k)| \le 4r|u-u_k|\le 4\sqrt{\frac{r\log n}{n}}.
\]
Now take any $u \ge u_N$. Since $m_i \in \{0,\pm 1,\ldots, \pm r\}$, therefore $|\tanh(um_i) - \tanh(u_N m_i)|\le 1-\tanh(u_N|m_i|)\le 1-\tanh(u_N)$. Thus,
\begin{align*}
|S(u)-S(u_N)| &\le \frac{4}{n}\sum_{i=1}^n |\tanh(um_i)-\tanh(u_N m_i)| \\
&\le 4(1-\tanh(u_N)) \le 4e^{-u_N} \le \frac{4e}{n}.
\end{align*}
If $n \ge 3$, then $\sqrt{\log n/n} \ge e/n$. Combining the steps, we see that for $n\ge 3$, 
\[
\min_{1\le k\le N} S(u_k) \le \min_{u\ge 0} S(u) + 4\sqrt{\frac{r\log n}{n}}.
\]
Finally, combining this with \eqref{tail1}, we get
\begin{align*}
&\pp\biggl\{S(\beta) \ge \min_{u\ge 0} S(u) + 4\sqrt{\frac{r\log n}{n}} + t\biggr\}\\
 &\le \pp\bigl\{S(\beta) \ge \min_{1\le k\le N} S(u_k) + t\bigr\}\\
 &\le \sum_{k=1}^N \pp\bigl\{S(\beta) \ge S(u_k) + t\bigr\} \le N \exp\biggl(-\frac{nt^2}{128(r+1)}\biggr).
\end{align*}
It is now easy to complete the proof by substituting the value of $N$ and choosing $t > \sqrt{Cr\log n/n}$ for sufficiently large $C$, so that the effect of $N$ washes out. 
\hfill $\Box$
\vskip.2in
\noindent Finally, let us prove our main result. 
\vskip.1in
\noindent {\bf Proof of Theorem \ref{conc}.} Let us begin with a useful general identity. Suppose $h:\xx \ra \rr$ is any measurable map such that $\ee|h(X) F(X,\xp)|<\infty$. Then clearly $\ee(h(X)f(X)) = \ee(h(X)F(X,\xp))$. 
Using the exchangeability of $X$ and $\xp$, and the antisymmetric
nature of $F$, we have
\[
\ee(h(X)F(X,\xp)) = \ee(h(\xp)F(\xp,X)) = -\ee(h(\xp)F(X,\xp)).
\]
Thus, we have
\begin{equation}\label{basis}
\ee(h(X) f(X)) = \ee(h(X)F(X,\xp)) = \frac{1}{2}\ee((h(X)-h(\xp))F(X,\xp)).
\end{equation}
The above equation is the basis of all that follows. First, note that by putting $h\equiv 1$, we immediately get $\ee(f(X))=0$, Similarly, part ($i$) of the Theorem follows by putting $h=f$. 
Next, let us start proving ($ii$). 
Let $m(\theta) := \ee(e^{\theta f(X)})$ be the moment generating function of $f(X)$. We can differentiate $m(\theta)$ and move the derivative inside the expectation because of the assumption that $\ee(e^{\theta f(X)} |F(X,\xp)|) < \infty$ for all~$\theta$. 
Thus, by equation~(\ref{basis}), we have
\begin{align*}
m^\prime(\theta) &= \ee(e^{\theta f(X)} f(X)) = \frac{1}{2} \ee((e^{\theta f(X)} - e^{\theta f(\xp)})F(X,\xp)).
\end{align*}
Now note that for any $x,y\in \rr$, 
\begin{equation}\label{expin}
\begin{split}
&\biggl|\frac{e^x-e^y}{x-y}\biggr| =\int_0^1 e^{tx + (1-t)y}dt \\
&\le \int_0^1 (te^x+(1-t)e^y)dt = \frac{1}{2}(e^x + e^y).
\end{split}
\end{equation}
Using this inequality, and the exchangeability of $X$ and $\xp$, we get
\begin{align*}
|m^\prime(\theta)| &\le \frac{|\theta|}{4}\ee((e^{\theta f(X)} + e^{\theta f(\xp)})|(f(X)-f(\xp)) F(X,\xp)|) \\
&= \frac{|\theta|}{2} \ee(e^{\theta f(X)} \Delta(X) + e^{\theta f(\xp)} \Delta(\xp))\\
&= |\theta|\ee(e^{\theta f(X)} \Delta(X))\\
&\le |\theta|\ee(e^{\theta f(X)} (Bf(X) + C)) = B|\theta| m^\prime(\theta) + C|\theta| m(\theta).
\end{align*}
Since $m$ is a convex function and $m^\prime(0) 
= \ee(f(X)) = 0$, therefore $m^\prime(\theta)$ always has the same sign as $\theta$. Thus, for $0\le \theta< 1/B$, the above inequality translates into 
\[
\frac{d}{d\theta} \log m(\theta) \le \frac{C\theta}{1-B\theta}.
\]
Using this and recalling that $m(0)=1$, we have
\begin{align*}
\log m(\theta) &\le \int_0^\theta \frac{Cu}{1-Bu}du \le \frac{C\theta^2}{2(1 - B\theta)}.
\end{align*}
Putting $\theta = t/(C+Bt)$, we get
\begin{align*}
\pp\{f(X) \ge t\} &\le \exp(-\theta t + \log m(\theta)) \le e^{-t^2/(2C + 2Bt)}.
\end{align*}
The lower tail can be done similarly; note that for $\theta \le 0$, we have $m^\prime(\theta)\le 0$, and hence
\[
|m^\prime(\theta)|\le B|\theta| m^\prime(\theta) + C|\theta| m(\theta)\le C|\theta| m(\theta),
\]
and this is the reason why $B$ does not appear in the lower tail bound. This completes the proof of part ($ii$). For the moment inequalities in part ($iii$), first observe that by equation (\ref{basis}), we have
\[
\ee(f(X)^{2k}) = \frac{1}{2} \ee((f(X)^{2k-1} - f(\xp)^{2k-1})F(X,\xp)).
\]
By the inequality 
\[
|x^{2k-1} - y^{2k-1}| \le \frac{2k-1}{2}(x^{2k-2} + y^{2k-2})|x-y|
\]
which follows easily from a convexity argument very similar to (\ref{expin}), we have
\begin{align*}
\ee(f(X)^{2k}) &\le (2k-1)\ee(f(X)^{2k-2} \Delta(X))
\end{align*}
By H\"{o}lder's inequality, we get
\[
\ee(f(X)^{2k}) \le (2k-1)(\ee(f(X)^{2k}))^{(k-1)/k}(\ee(\Delta(X)^k))^{1/k}.
\]
The proof is completed by transferring $\ee(f(X)^{2k})^{(k-1)/k}$ to the other side.\hfill $\Box$

\end{document}